\documentclass[12pt]{amsart}
\usepackage{amsmath,amsfonts,euscript,amscd,amsthm,amssymb,upref,graphics}

\theoremstyle{plain}
\newtheorem{theorem}{Theorem}
\swapnumbers

\newtheorem{lemma}[subsection]{Lemma}

\theoremstyle{definition}
\newtheorem{definition}[subsection]{Definition}

\newtheorem{remark}[subsection]{Remark}

\newtheorem{nothing*}[subsection]{}

\newcommand{\rien}[1]{}

\newcommand{\Aut}{ \operatorname{{\rm Aut}}}

\newcommand{\Nat}{\ensuremath{\mathbb{N}}}

\newcommand{\C}{\ensuremath{\mathbb{C}}}
\newcommand{\Reals}{\ensuremath{\mathbb{R}}}

\newcommand{\cO}{{\ensuremath{\mathcal{O}}}}

\def\e{\epsilon}

\newcommand{\emb}{\hookrightarrow}

\renewcommand{\epsilon}{\varepsilon}
\renewcommand{\phi}{\varphi}
\renewcommand{\emptyset}{\varnothing}
\begin{document}
\renewcommand{\baselinestretch}{1.07}

\title[Embedding some Riemann surfaces into $\C^2$ with interpolation ]
{Embedding some Riemann surfaces into $\C^2$ with interpolation}
\author{Frank Kutzschebauch}
\address{Mathematisches Institut \\ Universit\"at Bern
    \\Sidlerstr. 5
   \\ CH-3012 Bern, Switzerland}
\email{Frank.Kutzschebauch@math.unibe.ch}
\author{Erik L\o w}
\address{Matematisk Institutt\\  Universitet i Oslo\\ Blindern\\ Oslo}
\email{elow@math.uio.no}
\author{Erlend Forn\ae ss Wold}
\address{Matematisk Institutt\\  Universitet i Oslo\\ Blindern\\ Oslo}
\email{erlend@math.uio.no}

\date{January 31, 2007}
\thanks{Kutzschebauch supported by Schweizerische Nationalfonds grant 200021-107477/1}

\begin{abstract}
We formalize a technique for embedding Riemann surfaces properly
into $\C^2$, and we generalize all known embedding results to allow
interpolation on prescribed discrete sequences.
\end{abstract}
\maketitle \vfuzz=2pt

\vfuzz=2pt
\section{Introduction}

It is known that every Stein manifold of dimension $n>1$ admits a
proper holomorphic embedding in $\C^N$ with
$N=\left[\frac{3n}{2}\right] + 1$,and this $N$ is the smallest
possible due to an example of Forster \cite{Fs1}. The corresponding
embedding result with $N$ replaced by
$N'=\left[\frac{3n+1}{2}\right] + 1$ was announced by Eliashberg and
Gromov in 1970 \cite{GE} and proved in 1992 \cite{EG}. For even
values of $n\in\Nat$ we have $N=N'$ and hence the result of
Eliashberg and Gromov is the best possible. For odd $n$ we have
$N'=N+1$, and in this case the optimal result was obtained by
Sch\"urmann \cite{Sch}, also for Stein spaces with bounded embedding
dimension. A key ingredient in these results is the homotopy
principle for holomorphic sections of elliptic submersions over
Stein manifolds \cite{Gro}, \cite{FP2}.

Combining the known embedding results and the theory of holomorphic
automorphisms of $\C^N$, Forstneri\v{c}, Ivarsson, Prezelj and the
first author \cite{FIKP}  proved the above mentioned embedding
results with additional interpolation on discrete sequences.

In the case $n=1$ the above mentioned methods do not apply. It is
still an open problem whether every open Riemann surface embeds
properly into $\C^2$. Recently the third author achieved results
concerning that problem, and in this paper we prove the
corresponding results with interpolation on discrete sequences,
thus solving the second part of Problem 1.6 in \cite{FIKP}:

Let $X$  be a Riemann surface.  We say that $X$  \emph{embeds into
$\C^2$  with interpolation} if the following holds for all
discrete sequences  $\{a_j\}\subset X$ and $\{b_j\}\subset\C^2$
without repetition: There exists a proper holomorphic embedding
$f\colon X\emb\C^2$  with $f(a_j)=b_j$  for $j=1,2,...$

\begin{theorem}
\label{Main} If X is one of the following Riemann surfaces then X
embeds into $\C^2$ with interpolation:
\begin{itemize}
\item[(1)]  A finitely connected planar domain.
\item[(2)]  A finitely connected planar domain with a regularly
convergent sequence of points removed.
\item[(3)]  A domain in a torus with at most two complementary
components.
\item[(4)]  A finitely connected subset of a torus whose complementary
components do not reduce to points.
\item[(5)]  A Riemann surface whose double is hyperelliptic.
\item[(6)]  A smoothly bounded Riemann surface in $\C^2$  with a
single boundary component.
\end{itemize}
\end{theorem}
This list includes all instances of open Riemann surfaces we are
aware of admitting proper holomorphic embeddings into $\C^2$.  We
note that embeddings of hyperelliptic Riemann surfaces were
obtained by Forstneri\v{c}  and \v{C}erne in \cite{CF}.

Our method of proof follows the idea of the third author in
\cite{FW2} to embed $X$ as a Riemann surface in $\C^2$ with
unbounded boundary components and then construct a Fatou Bieberbach
domain whose intersection with the closure is exactly $X$. We
construct the Fatou-Bieberbach domain not as a basin of attraction,
but as the set where a certain sequence of holomorphic automorphisms
of $\C^2$ converges. The main ingredient is a version of Lemma 1 in
\cite{FW1}. We tried to formalize the ingredients in the proof and
formulate a more general  technical theorem (Theorem
\ref{formalized}) which implies Theorem \ref{Main}.

At each step of the inductive construction we take  care of the
additional interpolation condition in the same clever way as in
\cite{FIKP}.

The above theorem has already been proved in the special case that
$X$ is the unit disc by Globevnik in \cite{Glo}, and in the special
case that $X$ is an algebraic curve in $\C^2$  by Forstneri\v{c},
Ivarsson, Prezelj and the first author in \cite{FIKP}.

More results on embedding with interpolation can be found in \cite{K}.
%
%
%
%
%

\section{Proof of theorem \ref{Main}}

We shall use the theory of holomorphic automorphisms of $\C^N$.

Let $\pi_i : \C^2 \to \C$ denote the projection onto the ith
coordinate, $\Delta_R$ denotes the open disc of radius $R$ in $\C$,
$\overline \Delta_R$ its closure, $B_R$ is the ball of radius $R$ in
$\C^2$ and $\overline B_R$ its closure.

\begin{definition}
Given finitely many disjoint smooth real curves in $\C^2$ without
self intersection $\Gamma_i = \{\gamma_i (t) : t \in [0, \infty )$
or $ t \in (-\infty,\infty )\} \quad i=1,2, \ldots , m$, $\Gamma =
\bigcup_{i=1}^m \Gamma_i$, and a countable subset $E\subset \C^2
\setminus \Gamma$, which is discrete in $\C^2 \setminus \Gamma$. We
say that $(\Gamma, E)$ has the {\sl nice projection property} if
there is a holomorphic automorphism $\alpha \in \Aut_{hol} (\C^2) $
of $\C^2$ such that,if $\beta_i (t)=\alpha (\gamma_i (t))$,
$\Gamma^\prime = \alpha (\Gamma)$ and $E^\prime = \alpha (E)$, then
the following holds:

\begin{itemize}

\item[(1)] $\lim_{|t|\to \infty }|\pi_1 (\beta_i (t)))| =\infty \quad i=1, 2, \ldots, m$
\item[(2)] There is an $M \in\Reals$ such that for all $R\ge M$  $\C\setminus (\pi_1
(\Gamma^\prime)) \cup \overline\Delta_R)$ does not contain any
relatively compact connected components.
\item[(3)] The restriction of the projection $\pi_1$ to $\Gamma^\prime \cup E^\prime$ is a proper map into $\C$.
\end{itemize}

\end{definition}
\begin{lemma}\label{fact}
Given a polynomially convex compact set $M\subset \C^N$ and a finite
or countably infinite set $E\subset \C^N\setminus M$ such that
$M\cup E$ is compact. Then $M\cup E$ is polynomially convex.
\end{lemma}
\begin{proof}
Let $z\in \C^N$ be an arbitrary point in the complement of $M\cup
E$. Choose a polynomially convex compact neighborhood $\tilde M$ of
$M$ which contains $M$ in its interior but does not contain the
point $z$. Observe that $E \setminus \tilde M$ consists of finitely
many points.

Let $f ,g \in \cO (\C^N)$ be a holomorphic functions with

$$
f(z) = 1,  \quad \sup_{w\in \tilde M} |f(w)| <{1 \over 2}
$$

and

$$
g(z) = 1, \quad g(w) = 0 \ \forall w\in E\setminus \tilde M.
$$

Then $h = f^ng$ satisfies $1=h(z)>\sup_{w\in M\cup E} |h(w)|$ for
$n$ sufficiently big. Thus $z$ does not belong to the polynomially
convex hull of $M\cup E$.
\end{proof}

\begin{lemma}
\label{Wolds lemma} Given a polynomially convex compact set
$K\subset \C^2$, a finite set of points $c_1, c_2, \ldots, c_l \in
K$ and  a ball $B$ containing $K$ and a positive number $\e >0$.
Moreover given a finite number of real curves $\Gamma=
\bigcup_{i=1}^m \Gamma_i$ and a discrete subset $E\subset
\C^2\setminus \Gamma$ as in the definition above having the nice
projection property , such that $(\Gamma \cup E )\cap K = \emptyset
$. Then there is an automorphism $\psi \in \Aut_{hol} (\C^2)$ such
that

\begin{itemize}
\item[(a)]  $\sup_{z\in K} |\psi(z)-z| <\e$
\item[(b)]  $\psi (c_i) = c_i \quad i=1, 2, \ldots, l$
\item[(c)]  $\psi (\Gamma \cup E ) \subset \C^2 \setminus B$
\end{itemize}

\end{lemma}

\begin{proof}
To simplify notation, we will assume we already have applied
$\alpha$, i.e. that (1),(2) and (3) hold with $\beta_i$,
$\Gamma^\prime$, $E^\prime$ replaced by $\gamma_i$, $\Gamma $, $E$.
It is clear that the result will follow by conjugating with $\alpha
$, if we choose a slightly larger polynomially convex set and a
sufficiently big ball.

In order to construct $\psi$ assume that $B = B_R$ where $R$ is so
big that $R>M$ (remember $M$ is from the nice projection property)
and $K \subset \overline\Delta_R \times \C$.

Set $\tilde \Gamma = \Gamma \cap ( \overline\Delta_R \times \C)$ and
$\tilde E = E \cap ( \overline\Delta_R \times \C)$. By the nice
projection property (3) $\tilde \Gamma \cup \tilde E$ is compact.
Take an isotopy of diffeomorphisms removing $\tilde \Gamma \cup
\tilde E$ from $\overline\triangle_R\times\C$ not intersecting $K$
at any time (first do it for the curves $\Gamma$, this will
automatically remove all points from $E$ except finitely many,  then
remove the finite number of remaining  points) and apply the
Andersen-Lempert theorem to $K \cup \tilde \Gamma \cup \tilde E$. By
a theorem of Stolzenberg \cite{Sb} $K \cup \tilde \Gamma $ is
polynomially convex. By Lemma \ref{fact} we have that $K \cup \tilde
\Gamma \cup \tilde E$ is also polynomially convex and the same is
true for all isotopies of that set.
%
%
We get a holomorphic automorphism $\phi \in \Aut_{hol} (\C^2)$ with

\begin{itemize}
\item[(a')]  $\sup_{z\in K} |\phi(z)-z| <{\e \over 2}$,
\item[(b')]  $\phi (c_i) = c_i, \quad i=1, 2, \ldots, l$,
\item[(c')]  $\phi (\tilde\Gamma \cup \tilde E ) \subset \C^2 \setminus (\overline\triangle_R\cup\C)$.
\end{itemize}

To achieve (b') either correct the approximating automorphism by
suitable shears which are small on a certain ball containing the
whole situation or working in the proof of the Andersen- Lempert
theorem with the geometric structure of vector fields vanishing on a
finite number of points in $\C^2$.

We will correct this automorphism $\phi$, which could move points
from $\Gamma \cup E $ into $B_R$ which have not been there before,
by precomposing it with a shear. For that we set $\Gamma E_R = \{ z
\in \Gamma \cup E : \phi (z) \in B= B_R\}$. By assumption the
complement of $\overline\Delta_R \cup \pi_1 (\Gamma \cup E)$ does
not contain any bounded component and by construction $\pi_1 (\Gamma
E_R) \subset \pi_1 (\Gamma \cup E)\setminus \Delta_R$.

%
%

Now define a Mergelyan setting like in the proof of Lemma 1 in
\cite{FW1} on $\overline\Delta_R \cup \pi_1 (\Gamma \cup E)$ to
construct a shear automorphism $s$ of $\C^2$ of the form $s(z, w)=
(z, w + f(z))$ which removes $\Gamma E_R$ from the compact
$\phi^{-1} (B_R)$ not bringing new points from the set $E$ into
$\phi^{-1} (B_R)$. To achieve this last property  the crucial facts
are first that the limit set of the  sequence $E$ is contained in
$\Gamma$ and second that the projection $\pi_1$ restricted to
$\Gamma \cup E$ is proper. Finally observe that the approximating
function $f$ from Mergelyan's theorem can be chosen to be zero at
the finite number of points  $\pi_1 (c_i) \quad i=1, 2, \ldots, l$
contained in $\Delta_R$ (where $f$ approximates zero). Finally, let
$\psi = \phi\circ s$.
\end{proof}

We shall refer to $K, F=\{c_1,\cdots ,c_l \}, B, \Gamma, E $ as data for the lemma.

\begin{definition}
An open Riemann surface $X \subset \C^2$ together with a discrete
sequence without repetition $A=\{a_j\}\subset X$ are called {\sl
suitable} if $X$ is a bordered submanifold of $\C^2$ such that
$\partial X$ is a collection $\partial_1 ,\cdots ,\partial_m$ of
unbounded smooth curves, and $(\Gamma ,A)$ satisfies the nice
projection property, where $\Gamma = \bigcup_{i=1}^m \partial_i$.
\end{definition}
%
%
%
%

Here is a lemma on polynomial convexity that is needed for our
main lemma:

\begin{lemma}
\label{pclemma} Let $X \subset \C^2$ be a bordered Riemann surface
with unbounded boundary components $\partial_1 ,\cdots ,\partial_m$.
Then there is an exhaustion $M_j$ of $X$ by polynomially convex
compact sets such that if $K \subset \C^2 \setminus \partial X$ is
polynomially convex and $K \cap X \subset M_i$ for some $i$, then $K
\cup M_i$ is polynomially convex.
\end{lemma}

\begin{proof}
This follows from (the proof of) Proposition 3.1 of \cite{FW2}.
\end{proof}

\begin{lemma}
\label{inductivestep} Let $X \subset \C^2$ be an open Riemann
surface and $A=\{a_j\}\subset X$ a discrete sequence without
repetition which are suitable. Let $\{b_j\} \subset \C^2$ be a
discrete sequence without repetition. Let $B\subset B'\subset \C^2$
be closed balls such that $\Gamma \cap B' = \emptyset$ and let $L=
X\cap B'$. Assume that if $b_j \in B\cup L$ then $b_j = a_j$ and if
$b_j \notin B\cup L$ then $a_j \notin B'$,i.e. $a_j \in X\setminus
L$. Given $\e>0$ and a compact set $K\subset X$, there exist a ball
$B''\subset\C^2$ containing $B'$ ($B''$ may be chosen as large as
desired), a compact polynomially convex set $M\subset X$ with $L\cup
K\subset M$, and a holomorphic automorphism $\theta$ of $\C^2$
satisfying the following properties:
\begin{itemize}
\item[(i)]    $|\theta(z)-z|<\e$ for all $z\in B\cup L$,
\item[(ii)]   If $a_j\in M$ for some index $j$ then $\theta(a_j)=b_j\in B''$,
\item[(iii)]  If $b_j\in B' \cup (\theta (X)\cap  B'') $, then $a_j \in
M$,
\item[(iv)]   $\theta(M) \subset {\rm Int} B''$,
\item[(v)]    If $a_j\in X\setminus M$ for some $j$ then $\theta(a_j)\notin
B''$,
\item[(vi)]   $\theta (\Gamma)\cap  B'' = \emptyset$.
\end{itemize}
\end{lemma}

\begin{remark}
\label{addendum} This is the fundamental inductive step of the
construction. Notice that the lemma states that the geometric
situation is preserved after applying $\theta$, i.e. that if $X, A,
B, B', \Gamma$ are replaced by $\theta (X), \theta (A), B', B'',
\theta (\Gamma )$, then the hypotheses of the lemma still hold.
\end{remark}

\begin{proof}
An automorphism $\theta$ with the required properties
will be constructed in two steps, $\theta=\psi\circ\phi$.

By Lemma \ref{pclemma} there is a polynomially convex compact set $M
\subset X$ such that $L \cup K \cup \{a_j ; b_j \in B' \} \subset M$
and $B \cup M$ is polynomially convex.

Since $B \cup L$ is also polynomially convex, by (repeated application of)
Proposition 2.1 of \cite{F1999} there is an automorphism $\phi$ such that

\begin{itemize}
\item[(a)]   $ |\phi(z)-z|< \frac{\e}{2}$ for all $z \in B \cup L$
\item[(b)]   $\phi(a_j)=b_j$ for all  $a_j\in M$
\end{itemize}

Now, $\phi (B \cup M)$ is polynomially convex and if $E' = \{ a_j ;
a_j \notin M \}$, then $(\phi (\Gamma),\phi (E'))$ has the nice
projection property. Let $B''$ be a large ball containing $\phi (B
\cup M) \cup B'$. By Lemma \ref{Wolds lemma}, applied to the data
$\phi (B \cup M) , F=\{ b_j ; a_j \in M \} , B'' , \phi (\Gamma),
\phi (E')$, there is an automorphism $\psi $ satisfying the
following:
\begin{itemize}
\item[(a')] $|\psi(w)-w|< \frac{\e}{2}$ when $w\in \phi(B \cup M)$,
\item[(b')] $\psi(b_j) = b_j$ for all $b_j\in F$
\item[(c')] $\psi(\phi (\Gamma) \cup \phi (E')) \subset \C^2 \setminus B''$
\end{itemize}

Let $\theta=\psi\circ\phi$. (i) follows from (a) and (a').
(ii) follows from (b) and (b'). (iv) follows from (a') and the definition of $B''$.
(v) follows from (c') and the definition of $E'$. (vi) follows from (c').

To prove (iii), notice that if $b_j \in B'$, then $a_j \in M$ by the
definition of $M$. Let $F' = \{ b_j \in B''\} \supset F$. If $b_j
\in \theta (X) \setminus F$, then there is a shear $s$ such that $s$
is close to the identity on $B''$, $s(b_i)=b_i$ for all $b_i \in F$
and such that $s\circ\theta (X)$ avoids $b_j$. Replacing $\theta$ by
$s\circ\theta $ does not destroy the other properties. Hence we may
assume that $\theta (X)$ avoids $b_j$ and therefore all points in
$F' \setminus F$. This implies (iii).

\end{proof}

\begin{theorem}
\label{formalized} Let $X$ be an open Riemann surface and
$A=\{a_j\}\subset X$ a discrete sequence without repetition which
are suitable. Let $\{b_j\} \subset \C^2$ be a discrete sequence
without repetition. Then there exists a proper holomorphic embedding
$f\colon X\emb \C^2$ satisfying $f(a_j)=b_j$ for $j=1,2,\ldots$.
\end{theorem}

 \begin{proof}
 Choose an exhaustion $K_1\subset K_2\subset \cdots \subset
\bigcup_{j=1}^\infty K_j=X$ by compact sets. Fix a number $\e$
with $0<\e<1$. We shall inductively construct the following:
\begin{itemize}
\item[(a)] a sequence of holomorphic automorphisms $\Phi_k$ of $\C^2$,
\item[(b)] an exhaustion $L_1\subset L_2\subset \cdots\subset \bigcup_{j=1}^\infty L_j=X$
by compact, polynomially convex sets (remember that a Runge subset
of $X$ is polynomially convex in $\C^2$),
\item[(c)] a sequence of balls
$B_1\subset B_2\subset\cdots \subset \bigcup_{j=1}^\infty B_j=\C^2$
centered at $0\in\C^2$ whose radii satisfy $r_{k+1}>r_k+1$
for $k=1,2,\ldots$,
\end{itemize}
such that the following hold for all $k=1,2,\ldots$
(conditions (iv) and (v) are vacuous for $k=1$):
\begin{itemize}
\item[(i)]  $\Phi_k(L_k)=\Phi_k(X) \cap B_{k+1}$,
\item[(ii)]  if $a_j\in L_k$ for some $j$ then $\Phi_k(a_j)=b_j$,
\item[(iii)] if $b_j \in \Phi_k(L_k)\cup B_k$ for some $j$ then $a_j \in L_k$, and so $\Phi_k(a_j)=b_j$,
\item[(iv)]  $L_{k-1}\cup K_{k-1} \subset {\rm Int}L_{k}$,
\item[(v)]  $|\Phi_{k}(z)-\Phi_{k-1}(z)|< \e\, 2^{-k}$ for all $z\in B_{k-1}\cup L_{k-1}$.
\item[(vi)] $\Phi_{k} (\Gamma) \subset \C^2\setminus B_{k+1}$
\end{itemize}

To begin we set $B_0=\emptyset$ and choose a pair of balls
$B_1\subset B_2 \subset \C^2$ whose radii satisfy $r_2\ge r_1+1$. A
combination of shears followed by an automorphism from Lemma
\ref{Wolds lemma} (the polynomial convex compact set consists here
of finitely many points only) furnishes an automorphism $\Phi_1$ of
$\C^2$ such that $\Phi_1(a_j)=b_j$ for all those (finitely many)
indices $j$ for which $b_j\in B_2$,  $\Phi_1(a_j)\in\C^2\setminus
B_2$ for the remaining indices $j$, and $\Phi_{1} (\Gamma) \subset
\C^2\setminus B_{2}$. Setting $L_1=\{z\in X\colon \Phi_1(z)\in
B_2\}$, the properties (i), (ii),  (iii)  and (vi) are satisfied for
$k=1$ and the remaining two properties (iv), (v) are void. Assume
inductively that we have already found sets $L_1,\ldots,L_k \subset
X$, that we have found balls $B_1,\ldots, B_{k+1}\subset\C^2$ and
automorphisms $\Phi_1,\ldots,\Phi_k$ such that (i)--(vi) hold up to
index $k$. We now apply Lemma \ref{inductivestep} with $B=B_k$,
$B'=B_{k+1}$, $X$ replaced by $X_k=\Phi_k(X) \subset\C^2$, $\Gamma$
replaced by $\Phi_k(\Gamma)$ and $L=\Phi_k(L_k) \subset X_k$. This
gives us a compact polynomially convex set $M=M_k\subset X_k$
containing $\Phi_k(K_k\cup L_k)$, an automorphism $\theta=\theta_k$
of $\C^2$, and a ball $B''=B_{k+2} \subset \C^2$ of radius $r_{k+2}
\ge r_{k+1}+1$ such that the conclusion of Lemma \ref{inductivestep}
holds. In particular, $\theta_k(M_k) \subset B_{k+2}$, the
interpolation condition is satisfied for all points $b_j \in
\theta_k(M_k)\cup B_{k+1}$, and the remaining points in the sequence
$\{\Phi_k(a_j)\}_{j\in\Nat}$ together with the curves $\Phi_k
(\Gamma)$ are sent by $\theta_k$ out of the ball $B_{k+2}$. Setting
$$
    \Phi_{k+1} =\theta_k\circ \Phi_k, \quad
        L_{k+1}= \{z\in X\colon \Phi_{k+1}(z)\in B_{k+2}\}
$$
one easily checks that the properties (i)--(vi) hold for the index
$k+1$ as well. The induction may now continue.

Let $\Omega\subset \C^2$ denote the set of points $z\in\C^2$ for
which the sequence $\{\Phi_k(z)\colon k\in \Nat\}$ remains bounded.
Proposition 5.2 in \cite{F1999} (p.\ 108) implies that
$\lim_{k\to\infty} \Phi_k=\Phi$  exists on $\Omega$, the convergence
is uniform on compacts in $\Omega$, and $\Phi\colon\Omega\to\C^2$ is
a biholomorphic map of $\Omega$ onto $\C^2$ (a Fatou-Bieberbach
map). In fact, $\Omega =\bigcup_{k=1}^\infty \Phi_k^{-1}(B_k)$
(Proposition 5.1 in \cite{F1999}). From (v) we see that $X\subset
\Omega$, from (vi) it follows that $\Gamma \cap \Omega = \emptyset$
i.e. $X$ is a closed subset of $\Omega$ implying that $\Phi$
restricted to $X$ gives a proper holomorphic embedding into $\C^2$.
Properties (ii), (iii) imply the interpolation condition
$\Phi(a_j)=b_j$ for all $j=1,2,\ldots$. This completes the proof of
the theorem.

\end{proof}

We may now prove the following generalization of Theorem 1 in
\cite{FW2}:
\begin{theorem}\label{gen}
Let $X\subset\C^2$  be a Riemann surface whose boundary components
are smooth Jordan curves $\partial_1,...,\partial_m$. Assume that
there are points $p_i\in\partial_i$  such that
$$
\pi_1^{-1}(\pi_1(p_i))\cap\overline X=p_i.
$$
Assume that $\overline X$  is a smoothly embedded surface, and
that all $p_i$  are regular points of the projection $\pi_1$. If
in addition $\pi_1(\overline X)\subset\C$  is bounded, then $X$
embeds into $\C^2$ with interpolation.
\end{theorem}
\begin{proof}
Let $\alpha_1,...,\alpha_m\in\C$  be constants and define the
following rational map $F:\C^2\rightarrow\C^2$:
$$
F(z,w)=(z,w+\sum_{j=1}^m\frac{\alpha_j}{z-p_j}).
$$
Let $\Gamma$  denote $\partial F(X)$.  It is not hard to see that
the constants $\alpha_j$ can be chosen such that $\Gamma$
satisfies the conditions on the curves in the definition of the
nice projection property (use the projection $\pi_2$).  Let
$A=\{a_j\}\subset X$ be a discrete sequence without repetition.
Since $\pi_1(X)$ is bounded it follows the $\pi_2$  is proper when
restricted to $\Gamma\cup A$ so the pair $(\Gamma,A)$  has the
nice projection property.  Thus $X$ and $\Gamma$  are suitable,
and the result follows from Theorem 2.
\end{proof}

 \smallskip
{\em Proof of theorem \ref{Main}.}  We start by proving (1): Let $X$
be a finitely connected planar domain, and let $z_1,...,z_k$ be the
complementary components of $X$  consisting of isolated points (if
such components exist).  Let $g\colon X\emb\C^2$  be the embedding
$$
g(z):=(z,\sum_{j=1}^k\frac{1}{z-z_i}).
$$
If there are no other complementary components than the points
$z_i$ then $(\partial g(X),A)$  has the nice projection projection
property for any discrete sequence $A=\{a_j\}\subset X$. If there
are more complementary components we may assume, by Koebe
uniformization, that $X$ is a circled subset of the unit disk, and
it is clear that $g(X)$ satisfies the condition in Theorem
\ref{gen}. \

To prove (2) let $S_1,...,S_m\subset\C$  be smooth compact slits
with an endpoint $q_i$  for each curve, let $L$ denote the closed
negative real axis, and let $\{z_j\}\subset\C$ be a discrete
sequence without repetition.  Let these sets be pairwise disjoint.
By uniformization results (see for instance \cite{gz}), we may
assume that $X$ is of one of the following two types
\begin{itemize}
\item[(a)]  $X=\C\setminus(\cup_{i=1}^m S_m\cup\{z_j\})$,
\item[(b)]  $X=\C\setminus(\cup_{i=1}^m S_m\cup\{z_j\}\cup L)$,
\end{itemize}
and that all $S_i$  are contained in $\triangle$. Let
$A=\{a_j\}\subset X$  be a discrete sequence without repetition. We
want to construct an embedding $f\colon X\emb\C^2$ such that the
boundary $\Gamma=\partial f(X)$ satisfies the conditions on the
curves in the definition of the nice projection property with
projection on the plane $z=w$, and such that $f(A)$ is contained in
the set
$$
D:=\{(z,w)\subset\C^2;|z|\leq 1 \ \mathrm{or} \ |w|\leq 1 \
\mathrm{or} \ |w|\geq 2|z|\}.
$$
In that case we see that the projection onto the plane $z=w$  is
proper when restricted to $\Gamma\cup f(A)$, i.e. that $f(X)$  and
$f(A)$ are suitable.  \

We will define $f$  as a mapping on the form
$f(\zeta)=(\zeta,h(\zeta)+g(\zeta))$ with
$h(\zeta)=\sum_{j=1}^m\frac{\alpha_j}{\zeta-q_j}$,
$g(\zeta)=\sum_{j=1}^\infty(\frac{\beta_j}{\zeta-z_j})^{N_j}$  for
some choice of constants $\alpha_j,\beta_j\in\C, N_j\in\mathbb{N}$.
\

We have that generic choices of coefficients $\alpha_j$ give that
the map $\zeta\mapsto(\zeta,h(\zeta))$  maps $X$  onto a surface
with a nice projection of the boundary curves onto the plane
$z=w$.  In particular we may choose them such that
$|h(a_j)|<\frac{1}{2}$ for all
$a_j\in\C\setminus\overline\triangle$.  If we next choose $g$ such
that $\|g\|_{\mathcal{C}^1(\partial X\setminus\{z_j\})}$ is small,
then $f$ will map $X$ onto a surface with a nice projection of the
boundary curves.  \

For each $j\in\mathbb{N}$  choose $\beta_j>0$  such that the disks
$\overline\triangle_j:=\overline\triangle_{\beta_j}(z_j)$  are
pairwise disjoint and such that $\overline\triangle_j\cap\partial
X\setminus\{z_j\}=\emptyset$ for all $j\in\mathbb{N}$.  Make sure
that $\partial\triangle_j\cap A=\emptyset$  for each $j$. Since
$(\frac{\beta_j}{\zeta-z_j})^{N}\rightarrow\infty$  as
$N\rightarrow\infty$  on $\triangle_j\setminus\{z_j\}$  and
$(\frac{\beta_j}{\zeta-z_j})^{N}\rightarrow 0$  as
$N\rightarrow\infty$  on $\C\setminus\overline\triangle_j$ it is
clear that we may choose the sequence $N_j$  such that $|g(a_i)|\geq
3|a_i|$  if $a_i\in\triangle_j$  for some $j$  and
$g(a_i)<\frac{1}{2}$  otherwise.  For any choice of $\delta>0$  we
may also choose the $N_j$'s such that $\|g\|_{\mathcal{C}^1(\partial
X\setminus\{z_j\})}<\delta$.  If $\delta$  is small enough then
$\partial f(X)$  has a nice projection onto the plane $z=w$  and
$f(A)\subset D$. \

To prove $(3)$ let $\lambda\in\C$  be contained in the upper half
plane and let $\varrho_\lambda$ be the Weierstrass $p-$function:
$$
\varrho_{\lambda}(\zeta):=\frac{1}{\zeta^2}+\sum_{m,n\in
\mathbb{N}^2\setminus\{0\}}\frac{1}{(\zeta-(m+n\lambda))^2}-\frac{1}{(m+n\lambda)^2}.
$$
If $2p\in\C$  is not contained in the lattice
$L_{\lambda}:=\{\zeta\in\C;\zeta=m+\lambda n\}$ we have that the
map
$$
\phi_p(\zeta):=(\varrho(\zeta),\varrho(\zeta-p)),
$$
determines is a proper holomorphic embedding of
$\mathbb{T}_{\lambda}\setminus\{[0],[p]\}$  into $\C^2$, where
$\mathbb{T}_{\lambda}$  is the torus obtained by dividing out $\C$
by the lattice group $L_{\lambda}$, and $[0]$  and $[p]$  are the
equivalence classes of the points $0$  and $p$ (see \cite{FW2} for
details). We treat three different cases: Assume first that the
complementary components of $X$ are two distinct points, i.e. $X$ is
some quotient $\mathbb{T}_{\lambda}$  with two points $[z_1]$ and
$[z_2]$ removed.  By a linear change of coordinates on $\C$ we may
assume that $z_2=-z_1$  and $2z_1\notin L_{\lambda}$.  Then
$\phi_{z_1}$ is a proper embedding of
$\mathbb{T}_{\lambda}\setminus\{[0],[z_1]\}$  into $\C^2$.  Now
$\varrho_{\lambda}(z_1)=\varrho_{\lambda}(z_2)=q\in\C$  so we may
chose a M\"{o}bius transformation $m:\mathbb{P}^1\emb\mathbb{P}^1$
such that $m(q)=\infty$.  We get that the map
$$
f(\zeta)=(f_1(\zeta),f_2(\zeta))=
(m\circ\varrho_{\lambda}(\zeta),\varrho_{\lambda}(\zeta-z_1))
$$
determines a proper embedding of
$X=\mathbb{T}\setminus\{[z_1],[z_2]\}$  into $\C^2$.  Moreover since
$\mathrm{lim}_{\zeta\rightarrow z_j} f_1(\zeta)=\infty$  for $j=1,2$
we have that the pair $(\partial f(X),f(A))$  has the nice
projection property for all discrete sequences $A=\{a_j\}\subset
X$.\

Next assume that the complementary components $K_1$  and $K_2$ of
$X$ are not both points.  If neither of them are points the result
follows from $(4)$  so we may assume that $K_1$ is the point $[0]$
and that $K_2$  is not a point. We may then assume that $K_2$ is a
smoothly bounded disk in $\mathbb{T}_{\lambda}$  and by choosing
$[p]\in K_2$ appropriately one sees that the map $\phi_p$ embeds $X$
onto a surface in $\C^2$  satisfying the conditions in Theorem
\ref{gen} above after the coordinate change $(z,w)\mapsto(w,z)$ (see
\cite{FW2}  for more details).

To prove $(4)$ we recall from Theorem 1' in \cite{FW3}  that a
subset of a torus without isolated points in the boundary embeds
onto a surface in $\C^2$  satisfying the conditions in Theorem
\ref{gen} above. \

Next let $X$  be as in $(5)$.  Then $X$  can be obtained by
removing a finite set $D_1,...,D_m$  of smoothly bounded
(topological) disks from a closed Riemann surface $\mathcal{R}$,
so $X=\mathcal{R}\setminus\cup_{i=1}^m\overline D_i$.  There
exists a separating pair of inner functions
$f,g\in\mathcal{A}(X)$, i.e. $f$  and $g$  separate points on
$\overline X$  and $|f(x)|=|g(x)|=1$  for all $x\in\partial X$
\cite{RU}\cite{GO}. Then the map $h:=(f,g)$  embeds $X$ properly
into the unit polydisk in $\C^2$, and by perturbing the boundary
$h(\partial X)$ slightly one obtains a surface satisfying the
conditions in Theorem \ref{gen}  above. \

If $X$ is a surface as in $(6)$ we have by the maximum principle
that either $X$  is a planar domain or the projection $\pi_1$
takes its maximum at a finite set of points
$q_1,...,q_s\in\partial X$. By a linear change of coordinates $X$
satisfies the conditions in Theorem \ref{gen}. \

\bibliographystyle{amsplain}

 \end{document}